\documentclass{bmcart}

\usepackage[utf8]{inputenc} 
\usepackage{amssymb}
\usepackage{amsthm}

\def\includegraphics{}

\usepackage{graphics}
\usepackage{epsfig}
\usepackage{graphicx}
\usepackage{epstopdf}

\usepackage{amsmath, amsthm, amscd, amsfonts, amssymb, graphicx, color}
\usepackage{amsthm}
\usepackage{amssymb}
\usepackage{multirow,multicol}
\usepackage{booktabs}
\usepackage{amsmath, amscd, amsfonts, amssymb, graphicx, color}

\usepackage{float}

\begin{document}

\begin{frontmatter}

\begin{fmbox}
\dochead{Research}


\title{Symplectic Euler scheme for Hamiltonian stochastic differential equations driven by Levy noise }


\author[
   email={zhan2017@fafu.edu.cn,qzhan3@iit.edu}   
]{\fnm{Qingyi} \snm{Zhan}}$^{1}$
\author[
]{Jinqiao Duan$^2$}
\author[
]{Xiaofan Li$^2$}



\begin{artnotes}

\end{artnotes}

\end{fmbox}


\begin{abstractbox}

\begin{abstract}
This paper proposes a general symplectic Euler scheme for a class of Hamiltonian stochastic differential equations driven by L$\acute{e}$vy noise in the sense of Marcus form. The convergence of the symplectic Euler scheme for this  Hamiltonian stochastic differential equations is investigated. Realizable numerical implementation of this scheme is also provided in details. Numerical experiments are presented to illustrate the effectiveness and superiority of the proposed method by the simulations of its orbits, symplectic structure and Hamlitonian.

\end{abstract}




\end{abstractbox}
%

\end{frontmatter}



\section{Introduction}

\ \ \ \ \ \ Nowadays the stochastic process driven by non-Gaussian noise has played an important role in the theory and application of random dynamics, which can be modeled by stochastic differential equations(SDEs) \cite{D. Applebaum, J.Duan}. Especially, some Hamiltonian SDEs driven by L$\acute{e}$vy noise in the sense of Marcus form, which can preserve the symplectic structure, have been paid more and more attention in numerical simulations of many natural phenomena, such as the long-time orbits of n-body problem of planets and motion of particles in a fluid \cite{K. Feng, E.Hairer, J. Hong}. Numerical computation is in the center in the investigations of dynamical behaviour of Hamiltonian SDEs. Therefore, we mainly investigate the reliability  and feasibility of numerical computations of Hamiltonian SDEs driven by L$\acute{e}$vy noise.

 This work is motivated by two facts. First, developing highly accurate numerical methods for SDEs driven by non-Gaussian noise continue to be an interesting topic. Symplectic numerical integration scheme about Hamiltonian SDEs driven by Gaussian noise are shown in \cite{G.Milstein}.
 Many useful contributions are made in developing special numerical methods and the corresponding numerical analysis of SDEs \cite{Milstein},
 \cite{T. Wang, X. Wang}, \cite{Q.Zhan}-\cite{Q.Zhan5}. Second,
 the construction of conditions which can preserve the Hamiltonian structure of SDEs driven by non-Gaussian noise has been presented in \cite{P. Wei}.  These results are the foundations of symplectic scheme of Hamiltonian SDEs in the sense of Marcus form.
 To the best of our knowledge, no investigations of numerical symplectic scheme of Hamiltonian SDEs in the sense of Marcus form exist in the literature until now.

In this work, we mainly focus on the design  and implementation of symplectic Euler scheme(SES) for Hamiltonian SDEs driven by L$\acute{e}$vy noise in the sense of Marcus form. We numerically compare the dynamical behaviors of SES with those using non-symplectic methods in three aspects: the Hamiltonian, the preservations of symplectic structure and the phase trajectory in a long time interval. All these are illustrated in the numerical experiments. For our purpose that the numerical experiments are realizable and  simply achieved by programming, the L$\acute{e}$vy noise is restricted to be compound Poisson noise with a special realization\cite{T. Li}.

The results in this work show that under certain appropriate assumptions the symplectic structure is almost preserved in the discrete case in presence of the discontinuous input of L$\acute{e}$vy noises, and the numerical solution from this scheme can simulate the dynamical behaviour of Hamiltonian SDEs more accurately than non-symplectic methods in the long time interval.

The rest of this paper is organized as follows. Section 2 deals with some preliminaries. In Section 3 the theoretical results of preservation of symplectic structure are summarized. Section 4 presents the details of the symplectic Euler scheme. Illustrative numerical experiments are included in Section 5, where we demonstrate that the numerical implementation methods can be applied to obtain the numerical simulations of the long-time orbits of Hamiltonian SDEs in the sense of Marcus form. Finally, Section 6 is to summarize the conclusions of the paper.

\section{Preliminaries}
\ \ \ \ \  Let $(\Omega, \mathcal{F}, \{\mathcal{F}_t\}_{t\geq 0},  \mathbb{P})$, be a filtered probability space.
 We assume that $L(t)$ is a $d$-dimensional square integrable L$\acute{e}$vy process with the generating triplet $(\gamma,A,\nu)$,
where $\gamma$ is a $d$-dimensional drift vector, $A$ is a symmetric non-negative definite $d \times d$ matrix, and $\nu $ is a radially symmetric L$\acute{e}$vy jump measure on $\mathbb{R}^d \backslash 0$.

We consider the Cauchy problem for Hamiltonian SDEs driven by non-Gauss L$\acute{e}$vy noises in the sense of Marcus form on $\mathbb{M}$ as follows,
\begin{equation}\label{2.1}
 dX(t)=V_0(X(t))dt+\sum_{r=1}^m V_r(X(t)) \diamond  dL^r(t),\ \ \  \ \ \ X(0):=X(t_0)=x\in \mathbb{M},
\end{equation}
where $X\in \mathbb{R}^d$,$V_r: \mathbb{R}^d \rightarrow\mathbb{R}^d,r=0,1,...,m $, is the Hamiltonian vector fields, and $\mathbb{M}$ is a smooth $d$-dimensional manifold.
Here the Marcus integral for SDEs(\ref{2.1}) through Marcus mapping is usually written as
$$X(t)=x+\int_0^t V_0(X(s))ds+\sum_{r=1}^m \int_0^t V_r(X(s-)) \diamond  dL^r(s),\ \ \  \ \ \ $$
 which is defined as
$$X(t)=x+\int_0^tV_0(X(s))ds$$
$$+\sum_{r=1}^m\int_0^tV_r(X(s-))\circ dL_c^r(s) +\sum_{r=1}^m\int_0^tV_r(X(s-))dL_d^r(s)  $$
$$+\sum_{r=1}^{m}\sum_{0\leq s\leq t}\Big[\Phi^r(\Delta L^r(s), V_r(X(s-)),X(s-))-X(s-)-V_r(X(s-))\Delta L^r(s)\Big],$$
where $0<t\leq T< +\infty$, $L_c(t)$ and $L_d(t)$ are the usual continuous and discontinuous parts of $L(t)$, that is, $L(t)=L_c(t)+L_d(t)$. The notation $\circ$  denotes the Stratonovitch differential. And the flow map $\Phi^r(l,v(x),x)$ is the value at $s=1$ of the solution defined through the ordinary differential equations
\begin{equation}
 \left\{
  \begin{array}{lcr}
   \frac{d \xi^r}{ds}=V_r(\xi^r)l, s\in[0,1],\\
   \\
    \xi^r(0)=x.
  \end{array} \right.
 \end{equation}

Let us write Hamiltonian SDEs of even dimension $d=2n$ in the form of
\begin{equation}\label{2.2}
\begin{split}
dP=-\frac{\partial H_0}{\partial Q}(P,Q)dt-\sum_{r=1}^m \frac{\partial H_r}{\partial Q}(P,Q)\diamond  dL^r(t),\ \ \  \ \ \ P(t_0)=p,\\
dQ=\frac{\partial H_0}{\partial P}(P,Q)dt+\sum_{r=1}^m \frac{\partial H_r}{\partial P}(P,Q)\diamond  dL^r(t),\ \ \  \ \ \ Q(t_0)=q,
 \end{split}
\end{equation}
where $X=(P,Q)$, $X_0=(p,q)$ with $\| X_0\|<+\infty$  and
$V_r=(-\frac{\partial H_r}{\partial Q}, \frac{\partial H_r}{\partial P})$, $r=0,1,2,...,m.$ Here the norm $\|\cdot\|$ is defined as $(\ref{2.3})$,
and $P,Q,p,q$ are $n$-dimensional column-vectors. We assume that the functions $V_r,r=0,1,2,...,m $, satisfy  the conditions in \cite{P. Wei} such that  Hamiltonian SDEs $(\ref{2.2})$ have a unique global solution, and the solution process is adapted and c$\grave{a}$dl$\grave{a}$g.

Throughout the rest of this paper, we use the following notations.

Let $\mathbb{L}^2(\Omega,\mathbb{P})$ be the space of all bounded square-integrable random variables $x:\Omega \rightarrow \mathbb{R}^d$. For any random vector $x=({x_1},{ x_2},...,{ x_d}) \in \mathbb{R}^d$, the norm of  $x$ is defined in the form of
\begin{equation}\label{2.3}
\|x\|_2=\Big[\int_\Omega[|x_1(\omega)|^2+|x_2(\omega)|^2+...+|x_d(\omega)|^2]d\mathbb{P}\Big]^{\frac{1}{2}}< \infty.
\end{equation}
For any stochastic process  $x(t,\omega) \in \mathbb{R}^d$, the norm of  $x(t,\omega)$ is defined as follows
$$\|x(t,\omega)\|_2=\sup_{t\in \mathbb{R}^+}\|x_t(\omega)\|_2< \infty.$$

In addition, we define the norm of random matrices as follows
\begin{equation}\label{2.4}
 \| G \|_{\mathbb{L}^2(\Omega,\mathbb{P})} =\Big[\mathbb{E}(|G|^2)\Big]^{\frac{1}{2}},
\end{equation}
where $G$ is a random matrix and $|\cdot|$ is the operator norm.

For simplicity in notations, the norms $\|\cdot\|_2$ and $ \| \cdot \|_{\mathbb{L}^2(\Omega,\mathbb{P})}$ are usually written as $\|\cdot\|$.

\section{Theoretical results on preservation of symplectic structure }

\subsection{Preservation of symplectic structure for Hamiltonian SDEs}
\ \ \ \ From \cite{P. Wei}  we have

{\bf Lemma 3.1}\emph{
The flow map $X_t$ of the Hamiltonian SDEs $(\ref{2.2})$ is symplectic, that is,
$$dP\wedge dQ=dp\wedge dq, i.e., \sum_{i=1}^ndP^i\wedge dQ^i= \sum_{i=1}^ndp^i\wedge dq^i,$$
where $dP\wedge dQ$ is a differential two-form, $P=(P^1,P^2,...,P^n),$ and $Q=(Q^1,Q^2,...,Q^n)$.}

\subsection{Preservation of symplectic structure for discretized Hamiltonian SDEs}
\ \ \ \  In this section we consider the Hamiltonian SDEs with additive noise in the form of,
\begin{equation}\label{2.5}
\begin{split}
 dP=-\sigma_0(P,Q)dt-\sum_{r=1}^m \sigma_r(t)\diamond  dL^r(t),\ \ \  \ \ \ P(t_0)=p,\\
 dQ=\gamma_0(P,Q)dt+\sum_{r=1}^m \gamma_r(t)\diamond  dL^r(t),\ \ \  \ \ \ Q(t_0)=q,
 \end{split}
\end{equation}
where $$\sigma_0(P,Q)= \frac{\partial H_0}{\partial Q}(P,Q),\gamma_0(P,Q)=\frac{\partial H_0}{\partial P}(P,Q),$$
$$\sigma_r(t)=\frac{\partial H_r}{\partial Q}(P,Q),\gamma_r(t)=\frac{\partial H_r}{\partial P}(P,Q),r=1,2,...,m.$$

We make the assumption as follows.

{\bf Assumption 1.}

\emph{The functions $\sigma_0$ and $\gamma_0$ satisfy the Lipschitz condition
$$|\sigma_0(X_1)-\sigma_0(X_2)|\leq K|X_1-X_2|, |\gamma_0(X_1)-\gamma_0(X_2)|\leq K|X_1-X_2|,$$
where $K$ is a constant,  and $X_i=(P_i,Q_i)\in \mathbb{M},i=1,2$.}

The exact solution $X_{t_j}:=(P_{t_j},Q_{t_j})$ of $(\ref{2.5})$ at the time $t_j$ is shown as

 \begin{equation}\label{2.6}
\begin{split}
 P_{t_{j+1}}=P_{t_j}- \int_{t_j}^{t_{j+1}} \sigma_0(P(s),Q(s))ds-\sum_{r=1}^m \int_{t_j}^{t_{j+1}} \sigma_r(s) \diamond  dL^r(s),\ P(t_0)=p,\\
 Q_{t_{j+1}}=Q_{t_j}+\int_{t_j}^{t_{j+1}} \gamma_0(P(s),Q(s))ds+\sum_{r=1}^m \int_{t_j}^{t_{j+1}} \gamma_r(s) \diamond  dL^r(s),\ Q(t_0)=q,
 \end{split}
\end{equation}
where the Marcus integral for SDEs(\ref{2.6}) is usually defined by
$$\int_{t_j}^{t_{j+1}}\sigma_r(s)\diamond  \Delta L^r(s)=\int_{t_j}^{t_{j+1}}\sigma_r(s)\Delta L^r(s)$$
$$+\sum_{t_j\leq s\leq t_{j+1}}\Big[\Phi_1^r(\Delta L^r(s),\sigma_r(s),P(s-))-P(s-)-\sigma_r(s)\Delta L^r(s)\Big],$$
and
$$\int_{t_j}^{t_{j+1}}\gamma_r(s)\diamond  \Delta L^r(s)=\int_{t_j}^{t_{j+1}}\gamma_r(s)\Delta L^r(s)$$
$$+\sum_{t_j\leq s\leq t_{j+1}}\Big[\Phi_2^r(\Delta L^r(s),\gamma_r(s),Q(s-))-Q(s-)-\gamma_r(s)\Delta L^r(s)\Big].$$

And the flow maps $\Phi_1^r(l,\sigma_r(s),P(s-))$ and $\Phi_2^r(l,\gamma_r(s),Q(s-))$ are the value at $\hat{s}=1$ of the solutions defined through the ordinary differential equations, respectively,
\begin{equation}
 \left\{
  \begin{array}{lcr}
   \frac{d \xi_1^r}{d\hat{s}}=\sigma_r(s)l, \xi_1^r(0)=P(t_j-), \hat{s}\in[0,1],\\
   \\
    \frac{d \xi_2^r}{d\hat{s}}=\gamma_r(s)l, \xi_2^r(0)=Q(t_j-), \hat{s}\in[0,1].
  \end{array} \right.
 \end{equation}

We construct a stochastic semi-implicit Euler scheme for $(\ref{2.5})$
 \begin{equation}\label{2.7}
\begin{split}
 P_{j+1}=P_j- \sigma_0(P_{j+1},Q_j) \Delta t_j-\sum_{r=1}^m \sigma_r(t_j)\diamond  \Delta L^r(t_j),\ P_0=P(t_0)=p,\\
 Q_{j+1}=Q_j+\gamma_0(P_{j+1},Q_j)\Delta t_j+\sum_{r=1}^m \gamma_r(t_j)\diamond  \Delta L^r(t_j),\ Q_0=Q(t_0)=q,
 \end{split}
\end{equation}
where $\Delta t_j=t_{j+1}-t_j$, $t_0<t_1<,...,<t_N$,   $\Delta L^r(t_j)=L^r(t_{j+1})-L^r(t_{j})$, and $j=0,1,...,N$. Here the Marcus integral for SDEs(\ref{2.7}) is usually defined by
$$\sigma_r(t_j)\diamond  \Delta L^r(t_j)=\sigma_r(t_j)\Delta L^r(t_j)$$
$$+\sum_{t_j\leq s\leq t_{j+1}}\Big[\Phi_1^r(\Delta L^r(s),\sigma_r(t_j),P_j)-P_j-\sigma_r(t_j)\Delta L^r(s)\Big],$$
and
$$\gamma_r(t_j)\diamond  \Delta L^r(t_j)=\gamma_r(t_j)\Delta L^r(t_j)$$
$$+\sum_{t_j\leq s\leq t_{j+1}}\Big[\Phi_2^r(\Delta L^r(s),\gamma_r(t_j),Q_j)-Q_j-\gamma_r(t_j)\Delta L^r(s)\Big].$$
And the flow maps $\Phi_1^r(l,\sigma_r(t_j),P(t_j))$ and $\Phi_2^r(l,\gamma_r(t_j),Q(t_j))$ are the value at $\hat{s}=1$ of the solutions defined through the ordinary differential equations, respectively,
\begin{equation}
 \left\{
  \begin{array}{lcr}
   \frac{d \xi_1^r}{d\hat{s}}=\sigma_r(t_j)l, \xi_1^r(0)=P_j, \hat{s}\in[0,1],\\
   \\
    \frac{d \xi_2^r}{d\hat{s}}=\gamma_r(t_j)l, \xi_2^r(0)=Q_j, \hat{s}\in[0,1].
  \end{array} \right.
 \end{equation}
That is, we freeze $\sigma(s)$, $\gamma(s)$, $P(s-)$ and $Q(s-)$ on the right hand sight as $\sigma(t_j)$, $\gamma(t_j)$, $P_j$ and $Q_j$, respectively, which is the idea of Euler discretization in ODEs. $\cite{T. Li}$

We are in the position of the theorem which will show that the scheme $(\ref{2.7})$ is symplectic.

{\bf Theorem 3.2}\emph{
The semi-implicit-Euler $(\ref{2.7})$ for the Hamiltonian SDEs with additive noise $(\ref{2.5})$ preserves symplectic structure.}

\begin{proof}

It follows from the definition of symplectic structure that we only need to prove that
$$dP_{j+1}\wedge dQ_{j+1}=dP_{j}\wedge dQ_{j},j=0,1,2,...,N,$$
where $dP_j$ and  $dQ_j$ are the differential of $P_j$ and $Q_j$, respectively.

 Take the differential with respect to $P$ of the first equation in SDEs $(\ref{2.7})$, we obtain that
 $$ dP_{j+1}=dP_j- \frac{\partial \sigma_0}{\partial P}(P_{j+1},Q_j)\Delta t_j dP_{j+1}- \frac{\partial \sigma_0}{\partial Q}(P_{j+1},Q_j) \Delta t_jdQ_{j},$$
 that is,
  $$\Big[\mathbb{I}+ \frac{\partial \sigma_0}{\partial P}(P_{j+1},Q_j)\Delta t_j \Big]dP_{j+1}=dP_j-\frac{\partial \sigma_0}{\partial Q}(P_{j+1},Q_j) \Delta t_jdQ_{j},$$
  where $\mathbb{I}$ is the $n \times n$ unit matrix.

  Similarly, from the second equation of $(\ref{2.7})$ we can obtain that
 $$dQ_{j+1}-\frac{\partial \gamma_0}{\partial P}(P_{j+1},Q_j) \Delta t_jdP_{j+1}= \Big[\mathbb{I}+ \frac{\partial \gamma_0}{\partial Q }(P_{j+1},Q_j)\Delta t_j \Big]dQ_{j}.$$
 Multiply the above two equations, and we have
  $$\Big[\mathbb{I}+ \frac{\partial \sigma_0}{\partial P}(P_{j+1},Q_j)\Delta t_j \Big]dP_{j+1}\wedge dQ_{j+1}=\Big[\mathbb{I}+ \frac{\partial \gamma_0}{\partial Q}(P_{j+1},Q_j)\Delta t_j \Big]dP_{j}\wedge dQ_{j}.$$
  By the definition of Hamiltonian SDEs $(\ref{2.5})$, we have
  $$\frac{\partial \sigma_0}{\partial P}(P_{j+1},Q_j)= \frac{\partial \gamma_0}{\partial Q}(P_{j+1},Q_j)=\frac{\partial^2 H_0}{\partial Q \partial P}(P_{j+1},Q_j),$$
 since the following inequality usually holds
$$\Big[\mathbb{I}+ \frac{\partial^2 H_0}{\partial Q \partial P}(P_{j+1},Q_j)\Delta t_j \Big]\neq 0.$$
Therefore, we have
$$dP_{j+1}\wedge dQ_{j+1}=dP_{j}\wedge dQ_{j},j=0,1,2,...,N.$$
This completes the proof of Theorem 3.2.

\end{proof}

{\bf Remark 3.3}\emph{
As we know, the phase flow of Hamiltonian SDEs $(\ref{2.2})$ preserves symplectic structure. It follows from the proof of Theorem 3.2 that any schemes, which can preserve the symplectic structure of the deterministic parts of Hamiltonian SDEs $(\ref{2.2})$, can  preserve the symplectic structure of Hamiltonian SDEs $(\ref{2.2})$. Due to the complexity of the numerical implementation, we only focus on the method $(\ref{2.7})$ in this paper. More general and high-order symplectic schemes for Hamiltonian SDEs driven by the additive L$\acute{e}$vy noises will be presented in our future work.}

\subsection{Convergence of symplectic Euler scheme}
{\bf Theorem 3.4}\emph{
If the inequality $1-8\sqrt{2}K^2\tau^2> 0$ holds, where $K$ is the Lipschitz constant in Assumption 1 and $$\tau=\max_j \Delta t_j,$$ the scheme $(\ref{2.7})$ for the Hamiltonian SDEs with additive noise $(\ref{2.5})$ based on one-step approximation is of the mean-square order of accuracy 1.}

\begin{proof}
Let $\hat{X}_j=(\hat{P}_j,\hat{Q}_j)$ be the explicit Euler approximation of the Hamiltonian SDEs $(\ref{2.5})$ at the time $t_j,j=0,1,2,...,N$, and we obtain
 \begin{equation}\label{2.8}
\begin{split}
\hat{P}_{j+1}=\hat{P}_j- \sigma_0(\hat{P}_{j},\hat{Q}_j) \Delta t_j-\sum_{r=1}^m \sigma_r(t_j)\diamond  \Delta L^r(t_j),\ \hat{P}_0=P(t_0)=p,\\
 \hat{Q}_{j+1}=\hat{Q}_j+\gamma_0(\hat{P}_{j},\hat{Q}_j)\Delta t_j+\sum_{r=1}^m \gamma_r(t_j)\diamond  \Delta L^r(t_j),\ \hat{Q}_0=Q(t_0)=q.
 \end{split}
\end{equation}
It follows from Theorem 3.3 in \cite{T. Li} that we have
\begin{equation}\label{2.9}
\begin{split}
\sup_{j\leq N}\|\hat{X}_j-X_{t_j}\|\leq C\tau,
 \end{split}
\end{equation}
where $C$ is a constant.

 Now we assume that $X_{j}=(P_{j},Q_{j})$ be the numerical solution of the Hamiltonian SDEs $(\ref{2.5})$ at the time $t_j,j=0,1,2,...,N$, which is obtained by the method $(\ref{2.7})$.
 And we define $$\hat{E}_j=X_j-\hat{X}_{j}, \ E_j=X_j-X_{t_j}.$$
Then it is clear that
$$\hat{E}_{j+1}=
X_{j+1}-\hat{X}_{j+1}=
\left(\begin{array}{c}
  P_{j+1}-\hat{P}_{j+1}\\
  Q_{j+1}-\hat{Q}_{j+1}
\end{array}\right)
$$
$$
=
\left(\begin{array}{c}
P_{j}-\hat{P}_{j}\\
  Q_{j}-\hat{Q}_{j}
\end{array}\right)
+ \Delta t_j
\left(\begin{array}{c}
  -\sigma_0(P_{j+1},Q_j)+\sigma_0(\hat{P}_{j},\hat{Q}_j)\\
  \gamma_0(P_{j+1},Q_j)-\gamma_0(\hat{P}_{j},\hat{Q}_j)
\end{array}\right).
$$
Using the Lipschitz condition, Assumption 1, we obtain
$$ \mathbb{E}|-\sigma_0(P_{j+1},Q_j)+\sigma_0(\hat{P}_{j},\hat{Q}_j)|^2
=\mathbb{E}|\sigma_0(P_{j+1},Q_j)-\sigma_0(\hat{P}_j,Q_j)+\sigma_0(\hat{P}_j,Q_j)-\sigma_0(\hat{P}_{j},\hat{Q}_j)|^2 $$

$$\leq 2\mathbb{E}\Big[|\sigma_0(P_{j+1},Q_j)-\sigma_0(\hat{P}_j,Q_j)|^2+|\sigma_0(\hat{P}_j,Q_j)-\sigma_0(\hat{P}_{j},\hat{Q}_j)|^2\Big]  $$

$$\leq 2K^2\mathbb{E}\Big[|P_{j+1}-\hat{P}_j|^2+|Q_j-\hat{Q}_j|^2\Big ] \leq 2K^2\Big[\mathbb{E}|P_{j+1}-\hat{P}_j|^2+\mathbb{E}|Q_j-\hat{Q}_j|^2\Big ]  $$

$$\leq 2K^2\Big [\mathbb{E}|P_{j+1}-\hat{P}_{j+1}+\hat{P}_{j+1}-\hat{P}_j|^2+\mathbb{E}|\hat{E}_j|^2\Big ] $$

$$\leq2 K^2\Big[2\mathbb{E}|P_{j+1}-\hat{P}_{j+1}|^2+2\mathbb{E}|\hat{P}_{j+1}-\hat{P}_j|^2+\mathbb{E}|\hat{E}_j|^2\Big ] $$

$$\leq 2K^2 \Big [2\mathbb{E}|\hat{E}_{j+1}|^2+C'+\mathbb{E}|\hat{E}_j|^2\Big ],$$
where according to Theorem 2.1 in \cite{A.FerreiroCastilla} it is obvious to obtain
$$2\mathbb{E}\Big[\sup_j|\hat{P}_{j+1}-\hat{P}_j|^2\Big]\leq C',$$ here $C'$ depends on $K$ and $T$ only.

By the same way, we obtain that

$$\mathbb{E}|\gamma_0(P_{j+1},Q_j)-\gamma_0(\hat{P}_{j},\hat{Q}_j)|^2\leq 2K^2\Big [2\mathbb{E}|\hat{E}_{j+1}|^2+C'+\mathbb{E}|\hat{E}_j|^2\Big ].$$
Therefore, due to the former conclusion, we have
$$\mathbb{E}|\hat{E}_{j+1}|^2\leq 2\mathbb{E}|\hat{E}_{j}|^2
+2\tau^2 \mathbb{E}
\left|\begin{array}{c}
  -\sigma_0(P_{j+1},Q_j)+\sigma_0(\hat{P}_{j},\hat{Q}_j)\\
  \gamma_0(P_{j+1},Q_j)-\gamma_0(\hat{P}_{j},\hat{Q}_j)
\end{array}\right|^2
$$
$$\leq 2\mathbb{E}|\hat{E}_{j}|^2 + 4\sqrt{2}K^2\tau^2\Big [2\mathbb{E}|\hat{E}_{j+1}|^2+C'+\mathbb{E}|\hat{E}_j|^2\Big ].$$
That is,
$$(1-8\sqrt{2}K^2\tau^2)\mathbb{E}|\hat{E}_{j+1}|^2\leq 2(1+2\sqrt{2}K^2\tau^2)\mathbb{E}|\hat{E}_{j}|^2+4\sqrt{2}K^2C'\tau^2.$$
It follows from the assumption and the discrete version of Gronwall lemma that we have
$$\mathbb{E}|\hat{E}_j|^2\leq C''\tau^2.$$

Therefore,  we have
$$\sup_{j\leq N}\|E_j\|=\sup_{j\leq N}\|X_j-X_{t_j}\|\leq \sup_{j\leq N}\Big[\|X_j-\hat{X}_j\|+\|\hat{X}_j-X_{t_j}\|\Big]\leq C\tau,$$
where the last inequality refers to $(\ref{2.9})$.

 This finishes the proof of Theorem 3.4.
\end{proof}

\section{Numerical implementation methods}

\subsection{Experiment setup}

  \ \ \ \ For the realizability and simplify in programming,  we let $W_t, t\in \mathbb{R}^+:=[0,+\infty)$ be a one dimensional Wiener process, and we assume that $L^r(t)\in \mathbb{R},r=1,2,...,m,$ is the underlying compound Poisson process and has the corresponding realization,
 \begin{equation}\label{4.00}
 L^r(t)=\sum_{k=1}^{N^r(t)}R_k^rH(t-\tau_k^r)+bW(t),r=1,2,...,m,
 \end{equation}
 where $\tau_k^r$ is the jump time with rate $\lambda$, $R_k^r\in \mathbb{R}$ is the jump size with distribution $\mu$ and $$\sup_k R_k^r<+\infty,$$ $N^r(t)$ is the number of jumps until time $t$, and
 $H(t)$ is the Heaviside function with unit jump at time zero.

 Therefore it follows from the the realization of L$\acute{e}$vy noises $(\ref{4.00})$ that the Marcus integral for SDEs(\ref{2.1}) through Marcus mapping is
 written as
$$X(t)=x+\int_0^tV_0(X(s))ds+b\sum_{r=1}^m\int_0^tV_r(X(s))\circ dW(s)$$
$$+\sum_{r=1}^{m}\sum_{k=1}^{N^r(t)}\Big[\Phi_g^r(X(\tau_k^r-),R_k^r)-X(\tau_k^r-)\Big],$$
where the flow map $\Phi_g^r$ at $t=\tau_k^r$ is defined through the ordinary differential equations
\begin{equation}\label{4.0}
 \left\{
  \begin{array}{lcr}
   \frac{d\xi^r}{ds}=V_r(\xi^r)R_k^r, s\in[0,1],\\
   \\
   \xi^r(0)=X(\tau_k^r-),\\
   \Phi_g^r(X(\tau_k^r-),R_k)=\xi^r(1).
  \end{array} \right.
 \end{equation}

With the above mathematical implementation of Marcus integral, we can simulate the orbits of Hamiltonian SDEs in the long time interval by symplectic Euler scheme(SES), which will be illustrated in detail in Section 4.2. Here we only consider the case $b=0$ in $(\ref{4.00})$, and we refer to the results which have been proposed in $\cite{T. Li, X. Wang}$. And the realization of the case $b\neq 0$ is much more complicated, which will be presented in our further work.

\subsection{Symplectic Euler scheme}
\ \ \  \  We denote $\exp(\lambda)$ as the exponentially distributed random variable with mean $\frac{1}{\lambda}$. And we present this algorithm for the Hamiltonian SDEs as follows,
\begin{equation}\label{4.1}
\begin{split}
 dP=-\sigma_0(P,Q)dt-\sum_{r=1}^m \sigma_r(t)\diamond  dL^r(t),\ \ \  \ \ \ P(t_0)=p,\\
 dQ=\gamma_0(P,Q)dt+\sum_{r=1}^m \gamma_r(t)\diamond  dL^r(t),\ \ \  \ \ \ Q(t_0)=q.
 \end{split}
\end{equation}

Step 1. Given $t=0$, initial value $(P_0,Q_0)$ and the end time $T$.

Step 2. Generate a waiting time $\tau\sim \exp(\lambda)$ and a jump size $R_r\sim \mu_r$, where $\mu_r(r=1,2,...,m)$ is the distribution of random jumps.

Step 3. Solve the following ODEs $(\ref{4.2})$ by symplectic Euler scheme with initial value $(P(t),Q(t))$ until time $s=\tau$ to get its solution $(P(u),Q(u)), u\in[t, t+\tau)$,
\begin{equation}\label{4.2}
\begin{split}
 dP=-\sigma_0(P,Q)dt,\ \ \  \ \ \ P(t_0)=p,\\
 dQ=\gamma_0(P,Q)dt,\ \ \  \ \ \ Q(t_0)=q.
 \end{split}
\end{equation}

Step 4. Solve the following ODEs $(\ref{4.3})$ with initial value $(P(t+\tau)-,Q(t+\tau)-)$ until time $s=1$ to get  $(P(t+\tau),Q(t+\tau))$,
\begin{equation}\label{4.3}
\begin{split}
&\frac{dx}{dt}=-\sum_{r=1}^m \sigma_r(t)R_r,\ \ \  \ \ \ x(0)=P((t+\tau)-),\\
& \frac{dy}{dt}=\sum_{r=1}^m \gamma_r(t)R_r,\ \ \  \ \ \ y(0)=Q((t+\tau)-).
 \end{split}
\end{equation}

Step 5. Set $t:=t+\tau$, and repeat Step 2 unless $t\geq T$.

\section{ Numerical experiments}

\ \ \ \ We consider the following SDEs  $\cite{P. Wei}$,i.e., linear stochastic oscillator with L$\acute{e}$vy noise,
\begin{equation}\label{5.1}
\begin{split}
 & dP=-Qdt+\beta \diamond dL_t, \ P(t_0)=0  ,\\
 & dQ=Pdt, \ Q(t_0)=1,
 \end{split}
\end{equation}
where $$H(P,Q)=\frac{1}{2}(P^2+Q^2), H_1(P,Q)=-\beta Q.$$
Obviously, it is a special linear Hamiltonian SDEs driven by additive L$\acute{e}$vy noise which is the same as we discussed.
For any given initial values $(P_0,Q_0)$, it follows from the results in $\cite{P. Wei}$ that the exact solution of SDEs $(\ref{5.1})$ is shown as
\begin{equation}\label{5.2}
\begin{split}
 & P(t)=P_0\cos t+Q_0\sin t+\beta \int_0^t\sin (t-s)dL_s,\\
 & Q(t)=-P_0\sin t+Q_0\cos t+\beta \int_0^t\cos (t-s)dL_s.
 \end{split}
\end{equation}
The SES of SDEs $(\ref{5.1})$ is written as
\begin{equation}\label{5.3}
\begin{split}
 & P_{j+1}=P_j-Q_jdt+\beta d L_j,\\
 & Q_{j+1}=Q_j+P_{j+1}dt.
 \end{split}
\end{equation}
The explicit Euler method(EEM) of SDEs $(\ref{5.1})$ is written as
\begin{equation}\label{5.31}
\begin{split}
 & P_{j+1}=P_j-Q_jdt+\beta d L_j,\\
 & Q_{j+1}=Q_j+P_jdt.
 \end{split}
\end{equation}

We focus on the dynamical behaviours of SES from three aspects: the preservation of symplectic structure, convergence of SES and the sepcial realization of the solution of SDEs $(\ref{5.1})$ by SES, which will be illustrated in Section 5.1-5.3. And in the realization of the L$\acute{e}$vy noise, we choose $L(t)$ to be a compound Poisson process with jump size which is simulated by the normal distribution $N(0,\sigma^2),\sigma=0.2$ and intensity $\lambda=5.0$.

\subsection{Preservation of symplectic structure of Hamiltonian SDEs $(\ref{5.1})$ }

\ \ \ \ The results of our numerical experiments are shown as Fig.1-3, which includes three parts: the comparison of sample trajectories, the evolution of domains in the phase plane and the conservation of the Hamiltonian obtained by EEM, SES and the exact solution.

To start we apply EEM and SES to Hamiltonian SDEs $(\ref{5.1})$, and we can compare the phase trajectories of numerical solutions obtained by EEM and SES with the exact solutions from $(\ref{5.2})$. In order to improve the accuracy of the comparison, the initial conditions  are the same, that is, the step size is $dt=0.08$, $T=20.0$,$\beta=1.0$, $N=500.0$ and the initial values is $P(0)=0,Q(0)=1.0$.

As we can see from Fig.1 that the approximations of a sample phase trajectory of Hamiltonian SDEs $(\ref{5.1})$ are simulated by the symplectic method $(\ref{5.3})$, SES, as well as the non-symplectic method, EEM, respectively. The exact phase trajectory $(\ref{5.2})$ is obtained, too.
 \begin{figure}[h]
   \centering
   \begin{minipage}{6.5cm}
       \includegraphics[width=3.8in, height=2.80in]{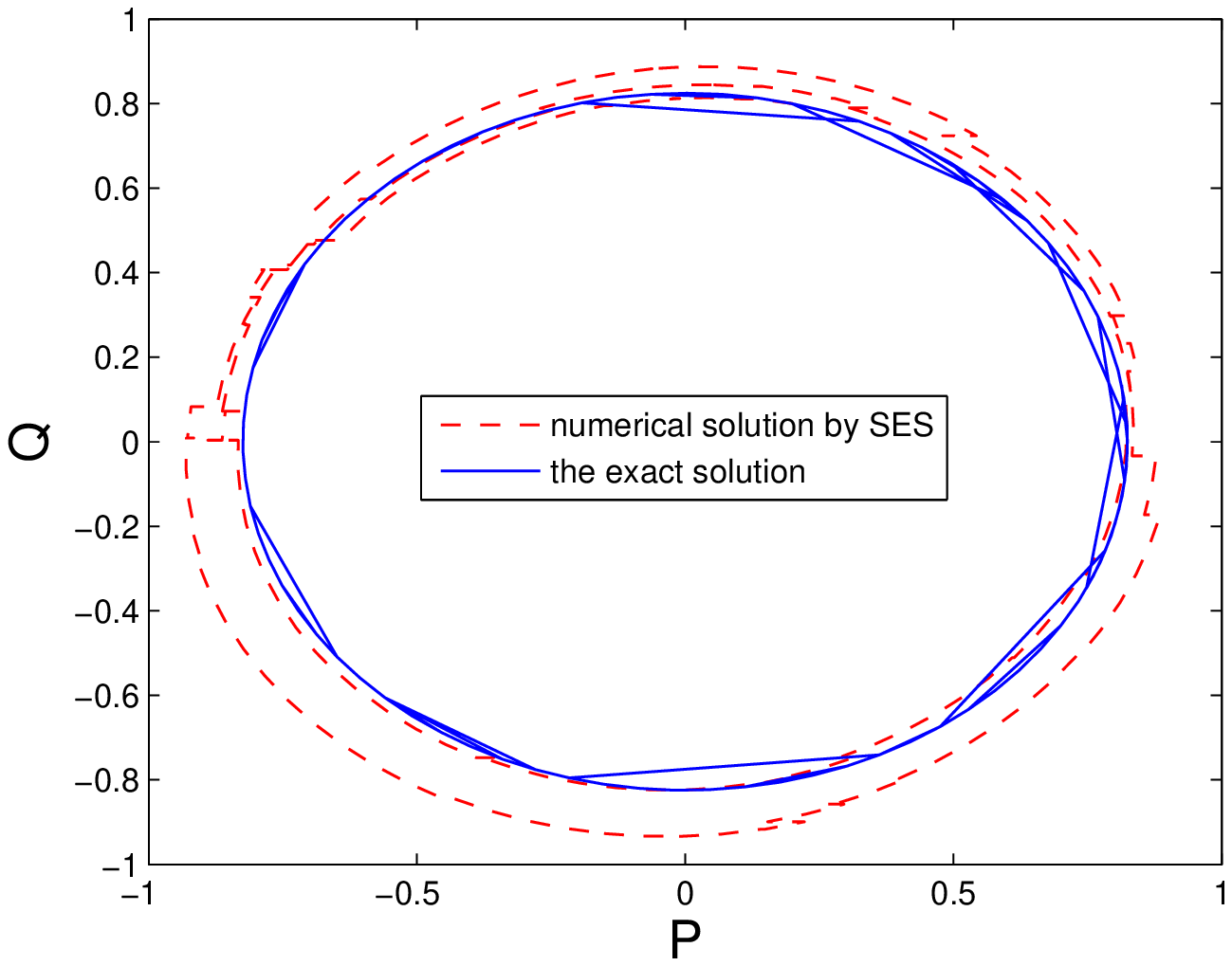}
    \end{minipage}
    \begin{minipage}{6.5cm}
       \includegraphics[width=3.8in, height=2.80in]{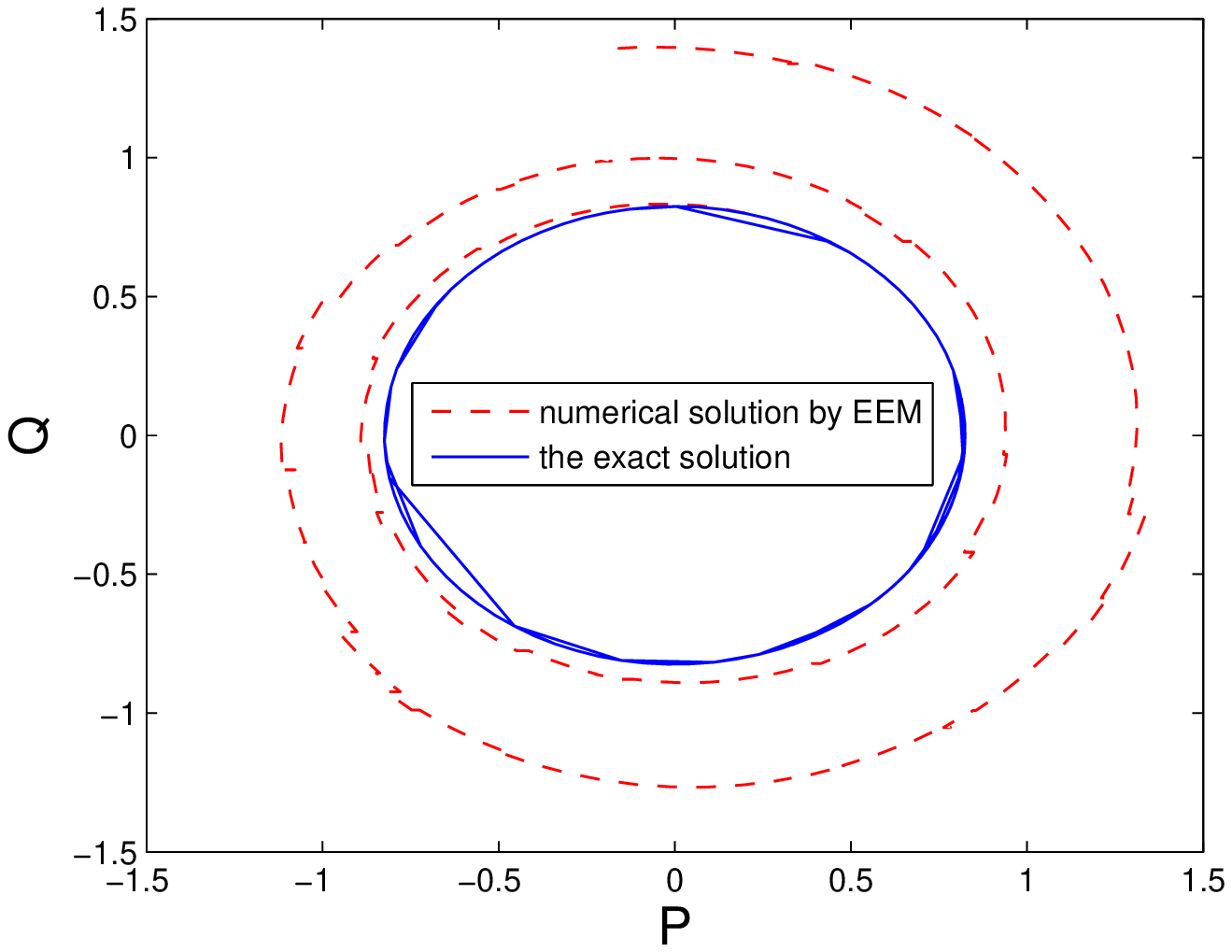}
    \end{minipage}
    \par{\scriptsize  Fig.1.  Comparison of zoom in parts of the exact trajectory of the solutions to SDEs $(\ref{5.1})$ obtained by $(\ref{5.2})$  and a sample trajectory obtained by SES(upper) and EEM(under), respectively.}
\end{figure}

 We find the fact that in the time interval $[0,20.0]$, the trajectory of the exact solution coincides almost well with that of SES, which is demonstrated in the upper panel of Fig.1, while the trajectory of EEM does not circle that of the exact solution, it disperses spirally  and quickly from the latter, which is shown as the under panel of Fig.1. It is obvious that SES has higher performance to preserve the circular phase trajectory than EEM. That is,  the structure of the trajectory of the solution to SDEs  $(\ref{5.1})$ obtained by EEM  obviously does not conserve the circular structure of that of the exact solution. The reason is that EEM has non-symplecticity, while SES dose.

This result indicates that EEM is unsuitable to simulate Hamiltonian SDEs $(\ref{5.1})$ in a long time interval. In contrast to EEM, SES reproduces the trajectory of SDEs $(\ref{5.1})$ more accurately.

Next we investigate the evolution of domains in the phase plane of SDEs $(\ref{5.1})$. Motivated by the work in Ref. $\cite{E.Hairer}$, we choose the
initial domain at the initial time $t=0$ with unit circle. In this section the initial conditions are revised as follows, the step size is $dt=0.08$, $T=20.0$,$\beta=1.0$, $N=500$ and the initial values is $P(0)=0.2,Q(0)=0.8$ for better comparison.

At three different time moments, $t=0$,$t=4.0$ and $t=8.0$, the images of these circles are demonstrated in the plane. These domains present the area of the phase space of points $(P,Q)$ at these time moments, and those points are on the trajectory obtained by the exact solution, SES and EEM, respectively.
 \begin{figure}[h]
   \centering
   \begin{minipage}{6.5cm}
       \includegraphics[width=3.8in, height=2.60in]{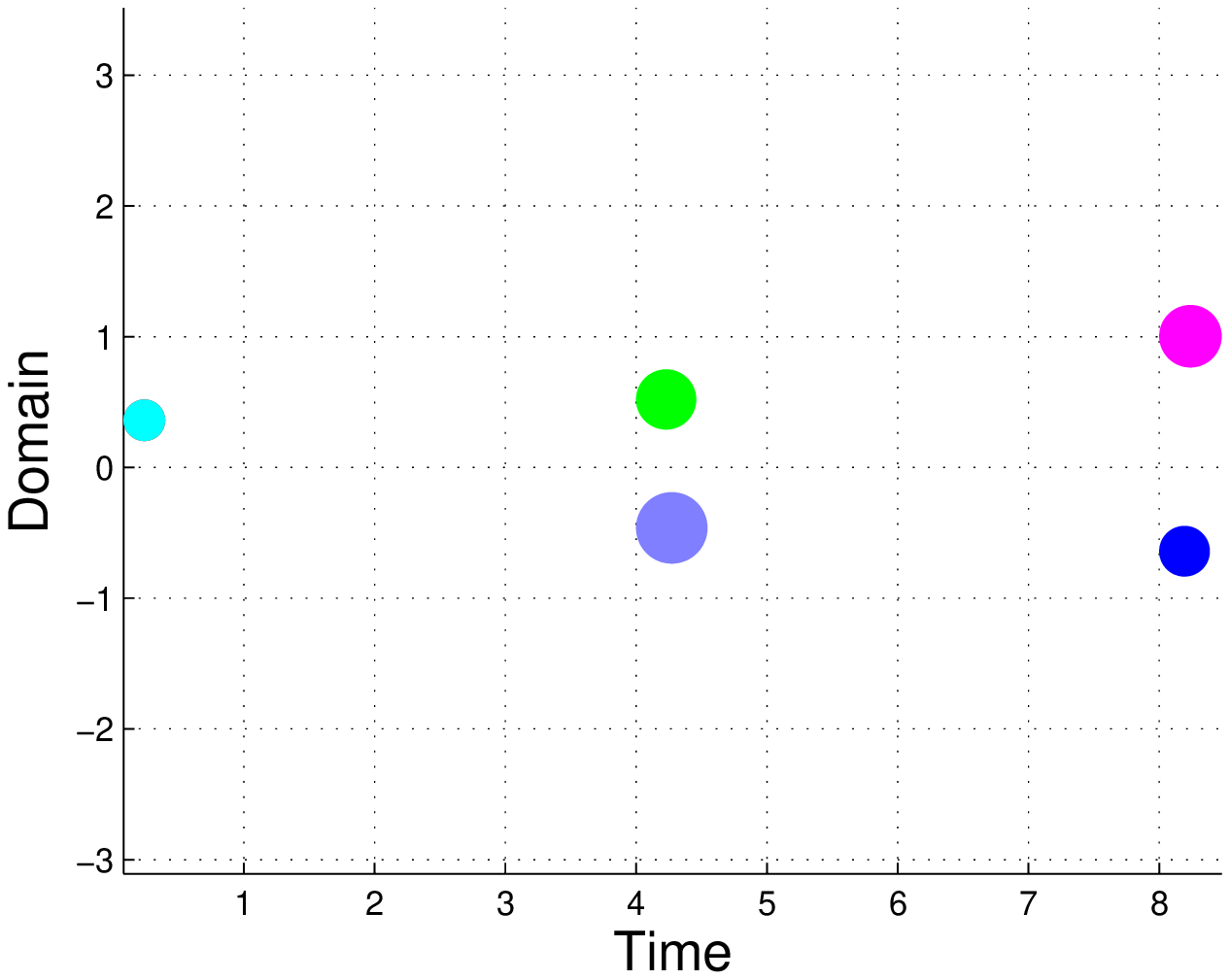}
    \end{minipage}
    \begin{minipage}{6.5cm}
       \includegraphics[width=3.8in, height=2.60in]{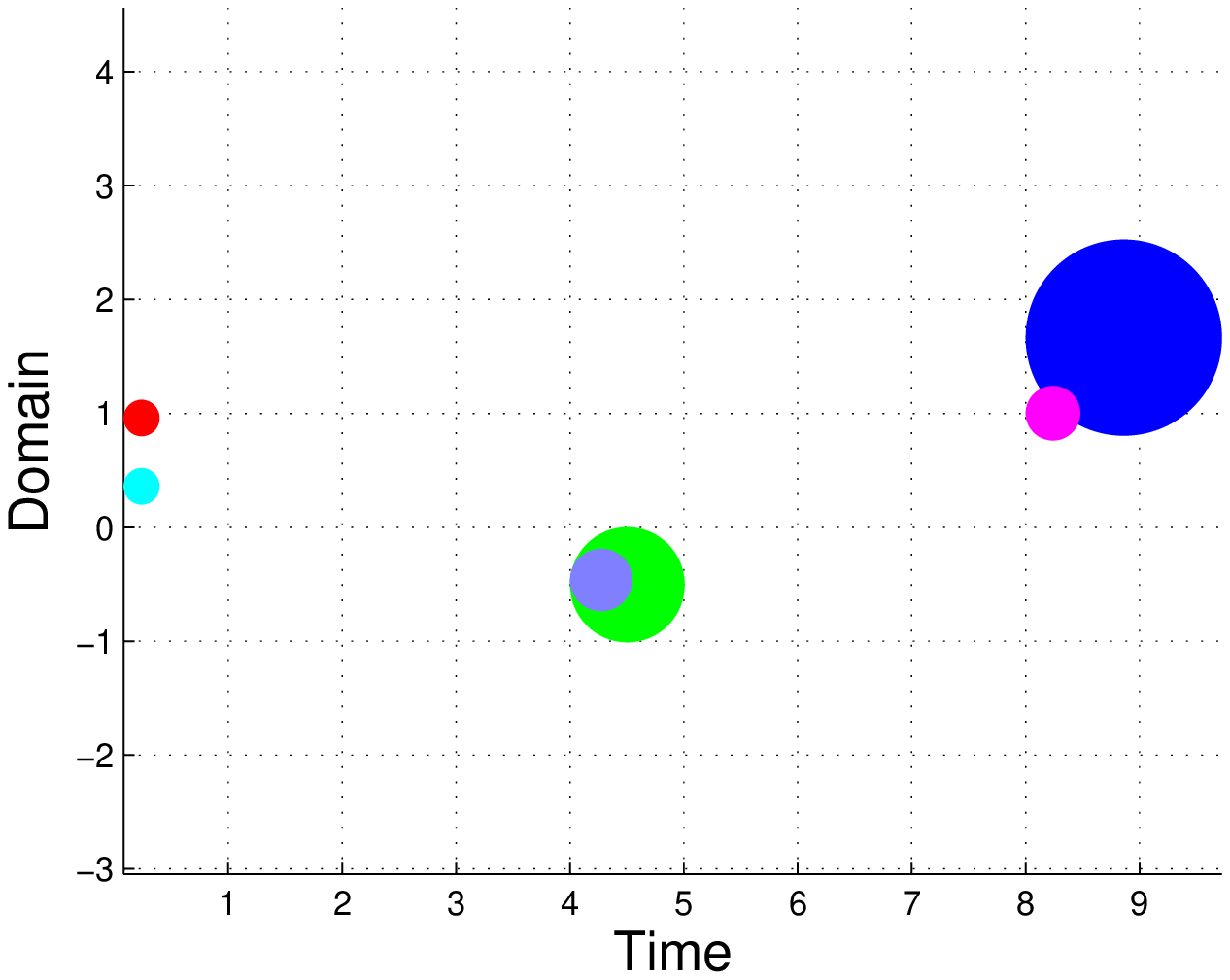}
    \end{minipage}
    \par{\scriptsize  Fig.2. Comparison of the domains in the phase plane of SDES $(\ref{5.1})$ obtained by the exact solution $(\ref{5.2})$, SES(upper) and EEM(under), respectively.}
\end{figure}

 As we can see, the images of the above three circles in the upper panel of Fig. 2 are obtained by the exact solution, and the ones of the below three circles are obtained by SES, where the image of the first circle is the same as the former. And it is clear that there is very little difference in the images between SES and the exact solution due to the influence of the L$\acute{e}$vy noise. This illustrates that the exact flow of a Hamiltonian SDEs $(\ref{5.1})$ can almost preserve the symplectic structure, which has been proved  theoretically in Theorem 3.3.

 On the contrary, in the under panel of Fig.2 the images of the above three circles  are obtained by EEM, and the ones of the below three circles are obtained by the exact solution. It is obvious that in the case of EEM, the images of these circles has the increasing radius such that the contrast on the images is significant. This fact is because of the reason that EEM dose not preserve symplectic structure. Despite the fact that EEM and SES have the same mean-square order of accuracy, SES has better performance on the approximate the exact images than EEM.

Lastly we check the Hamiltonians of SDEs $(\ref{5.1})$.
 \begin{figure}[h]
\centering
\includegraphics[width=3.8in, height=2.60in]{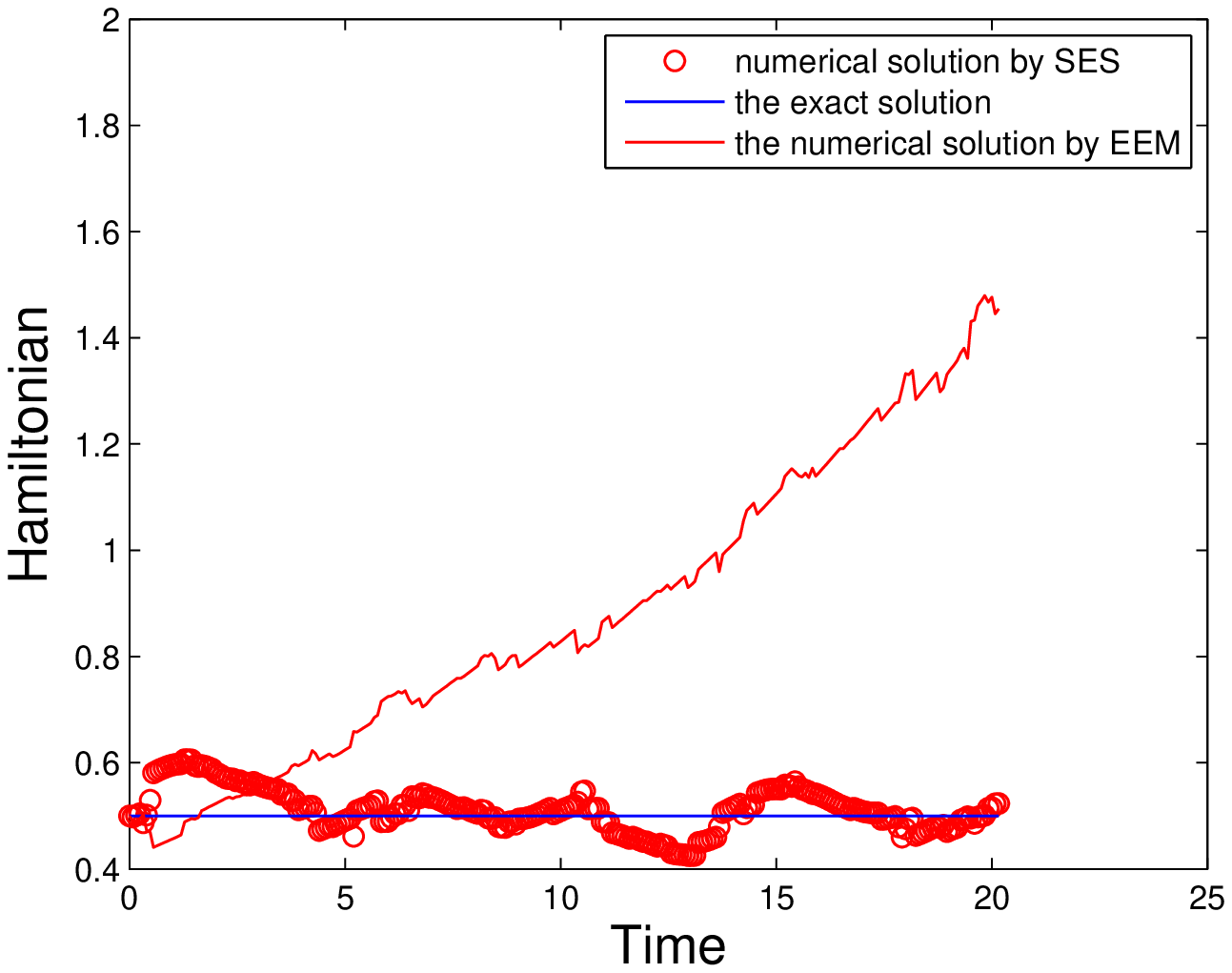}
\par{\scriptsize  Fig.3. Conservation of the Hamiltonian of SES, the exact solution and EEM.}
\end{figure}
It can be seen from Fig.3 that $H(P,Q)$ is an invariant of the exact solution of SDEs $(\ref{5.1})$. Due to the L$\acute{e}$vy noise, it can be approximately preserve by SES, that is, the curve of Hamiltonian jumps around the line Hamiltonian$=0.5$. However, non-symplectic numerical scheme, EEM dose not has this property such that the Hamiltonian increases indefinitely, which is shown as Fig.3. Here we take the initial conditions $T=20.0$, $P(0)=0.2$ and $Q(0)=0.8$.

\subsection{Convergence of SES}
\ \ \ \ This numerical experiment examines the convergence of SES. It is not difficult to see from Fig.4 that the convergence rate satisfies the inequality  log($\|$error of SES$\|$)$\leq 0.5$ for the end time $T=20.0$. Due to discontinuous inputting of the L$\acute{e}$vy noise, the curve has some jumps in some uncertain time moments, but it almost lays down the straight line log($\|$error of SES$\|$)$=0.5$. And  These phenomena verify the results of Theorem 3.5 that the mean-square order of the proposed method is 1. In this test we choose the same parameters as Section 5.1, the mean-square norm is taken as $(\ref{2.3})$.

\begin{figure}[h]
\centering
\includegraphics[width=3.8in, height=2.60in]{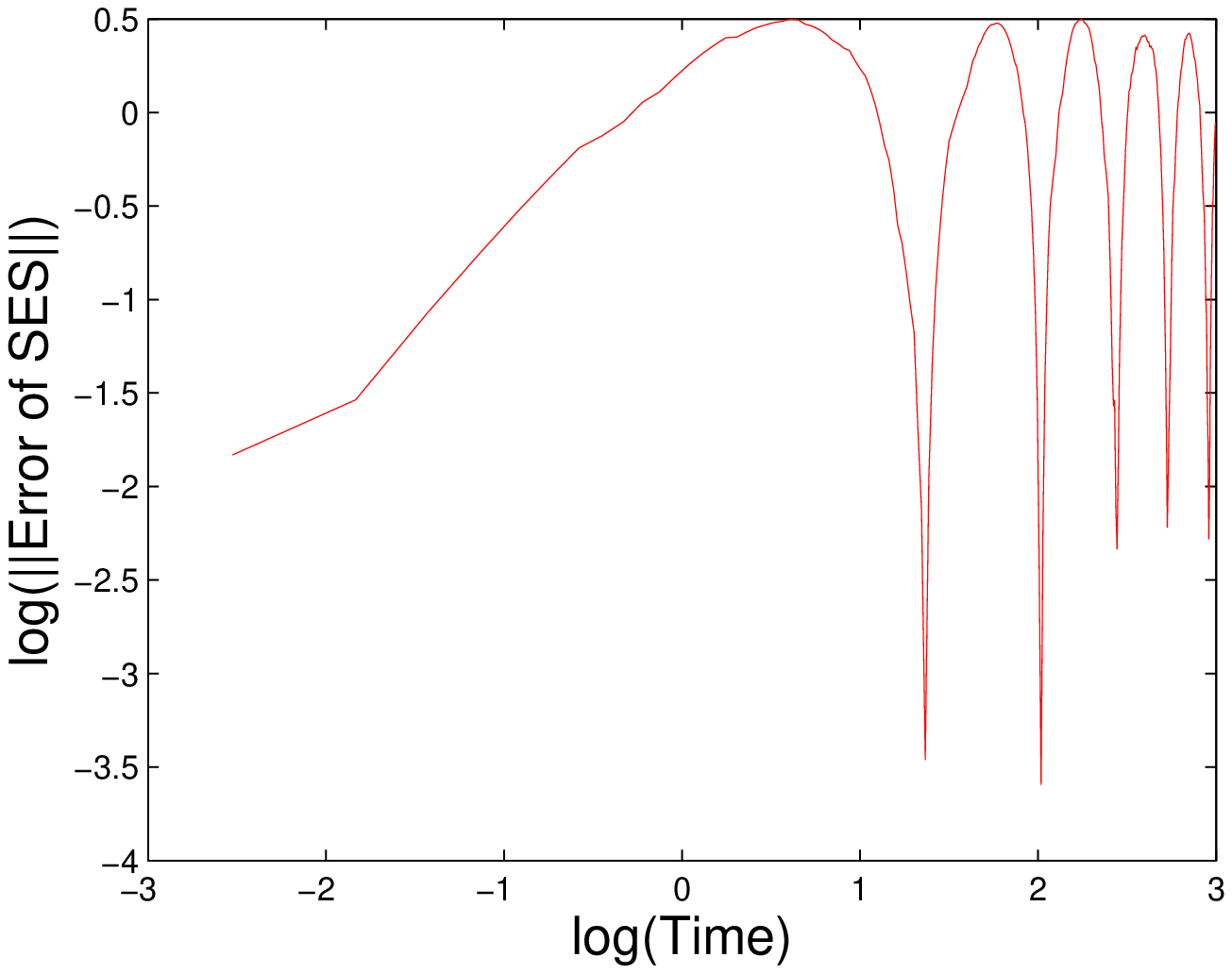}
\par{\scriptsize  Fig.4. The mean-square convergence rate of SES.}
\end{figure}

\subsection{Sepcial realization of the solution of SDEs $(\ref{5.1})$ by SES}

\ \ \ \ At first we will show the detail of a sample phase trajectory of SDEs $(\ref{5.1})$ simulated by SES, which can be viewed as the additional part of Section 5.1. In this experiment the red cycle symbol shows the numerical solution of SDEs $(\ref{5.1})$ obtained by SES. The upper panel in Fig.5 is  in the view of 3-dimensions, while for better comparing with the time $t$, the under panel in Fig.5 is  of 2-dimensions.

We can observe in Fig.5 that there are several discontinuous jumps in this numerical solution of SDEs $(\ref{5.1})$  in the interval $[0,1]$.  This verify the fact that SDEs $(\ref{5.1})$ is driven by L$\acute{e}$vy noise.
The initial conditions are similar to  which are presented in Section 5.1.
Some other parameters are $\Delta t=0.01$, $\beta=1.0$ and $T=1.0$.
 \begin{figure}[h]
   \centering
   \begin{minipage}{6.5cm}
       \includegraphics[width=3.8in, height=2.60in]{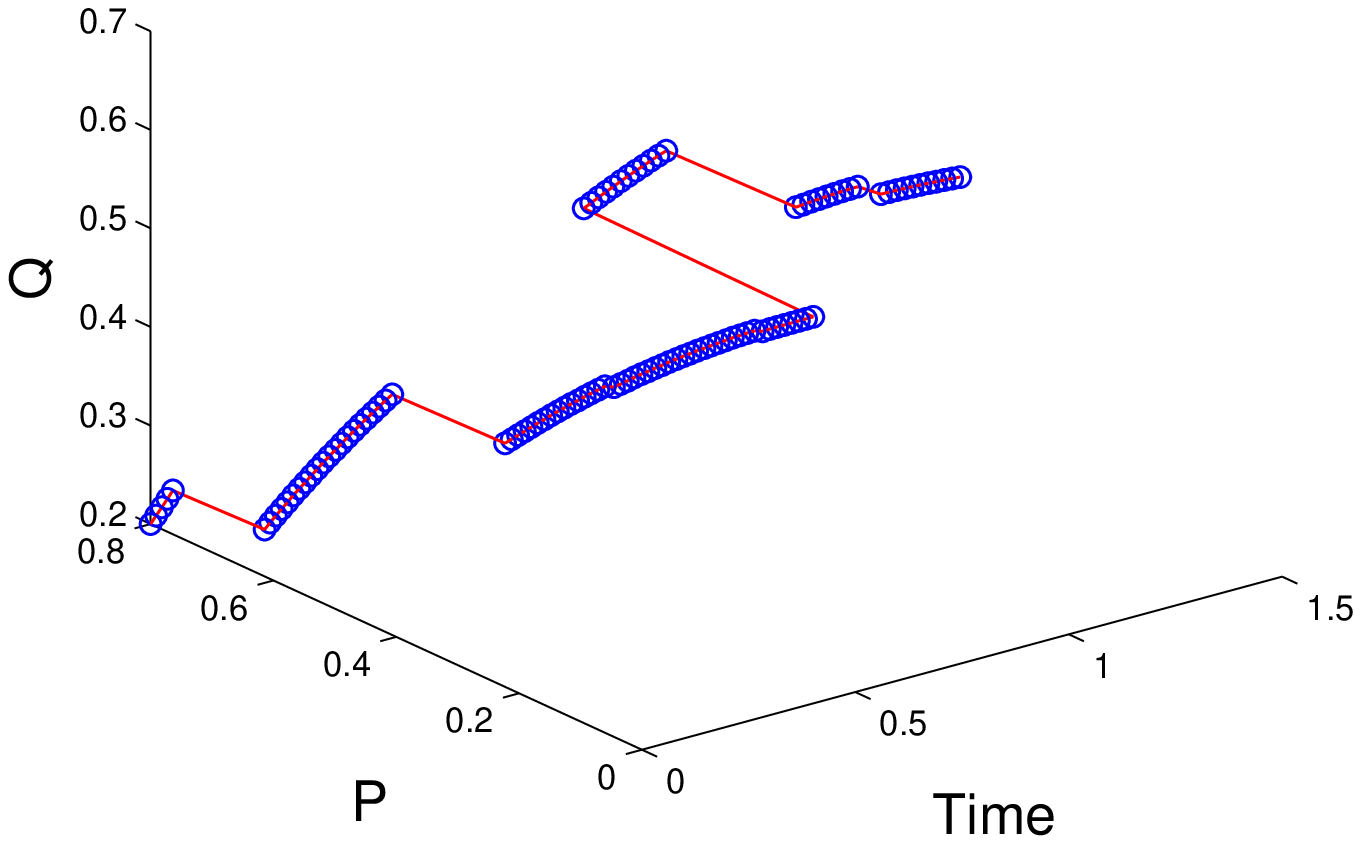}
    \end{minipage}
    \begin{minipage}{6.5cm}
       \includegraphics[width=3.8in, height=2.60in]{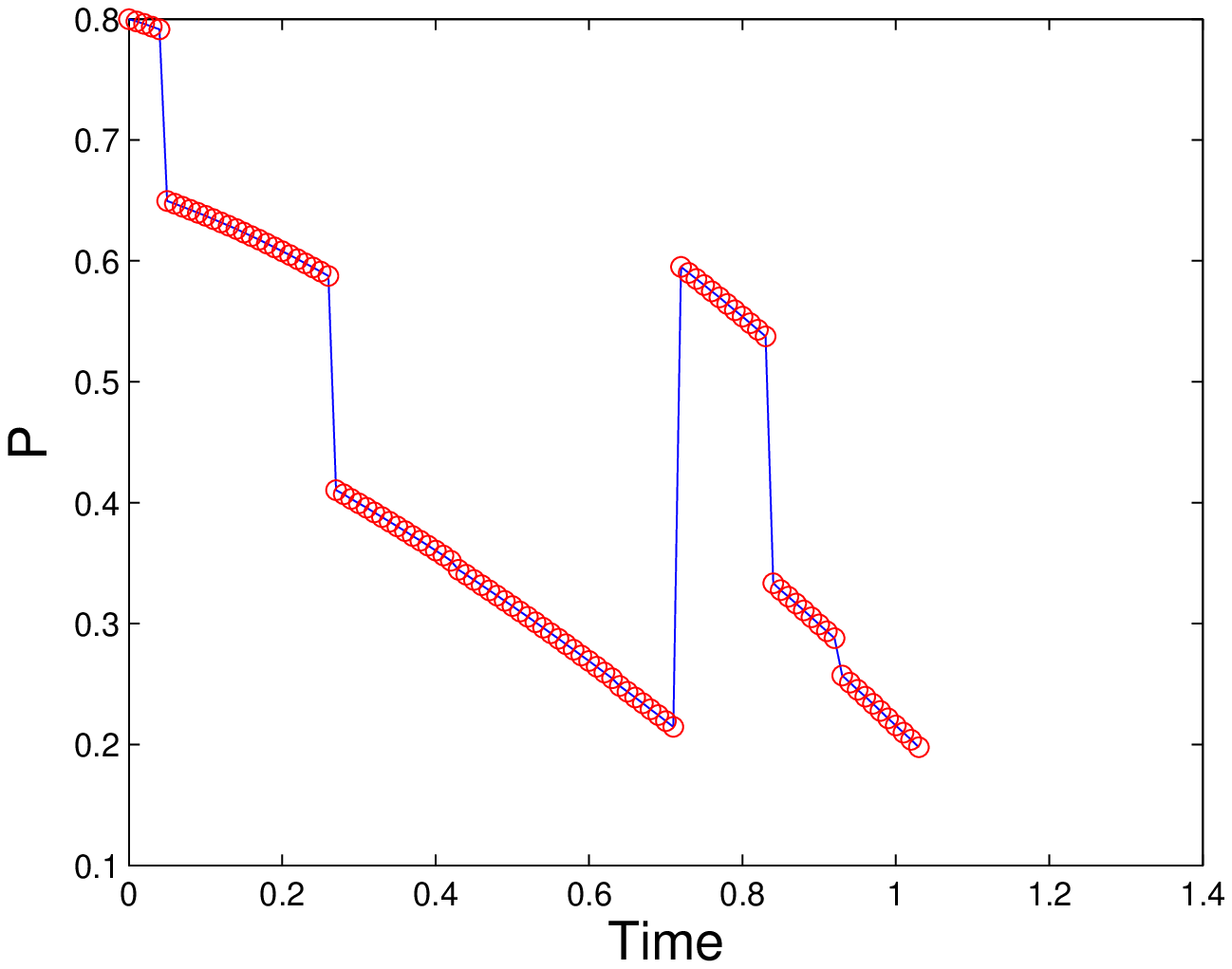}
    \end{minipage}
    \par{\scriptsize  Fig.5. Part zoom of one specific realization of the solution of SDEs $(\ref{5.1})$ by symplectic Euler method.}
\end{figure}

Next we will show the special realization of conservation of the Hamlitonian of the numerical solution of SDEs $(\ref{5.1})$ obtained by SES.
In this experiment the red cycle symbol presents the Hamiltonian of the  numerical solution of SDEs $(\ref{5.1})$ obtained by SES. And the under panel of Fig.6 presents the domain of three selected points in this numerical solution of SDEs $(\ref{5.1})$ obtained by SES at three time moments $t=0$, $t=0.5$ and $t=1.0$, respectively.

We can observe in the upper panel of Fig.6 that there are many jumps in this curve  in the interval $[0,1]$, which can explain the fact that the Hamiltonian of the numerical solution of SDEs $(\ref{5.1})$ by SES is not preserved as the one of the exact solution, but its mean is almost the same to the latter.  It is obvious in the under panel of Fig.6 that the images of the circles are almost the same at three different time moments. This demonstrates that the domians are almost the same in spite of some jumps.  In this experiment some parameters are
 $\Delta t=0.01$ and $T=1.0$.

 We summarize that these experiments demonstrate the better behaviour of the efficiency and superiority of SES than that of EEM, which has the same mean-square convergence order. Meanwhile, these also show the fact that symplectic methods are more suitable to compute the numerical solution of
 Hamiltonian SDEs with L$\acute{e}$vy noise in the sense of Marcus rather than non-symplectic methods.

 \begin{figure}[h]
   \centering
   \begin{minipage}{6.5cm}
       \includegraphics[width=3.8in, height=2.60in]{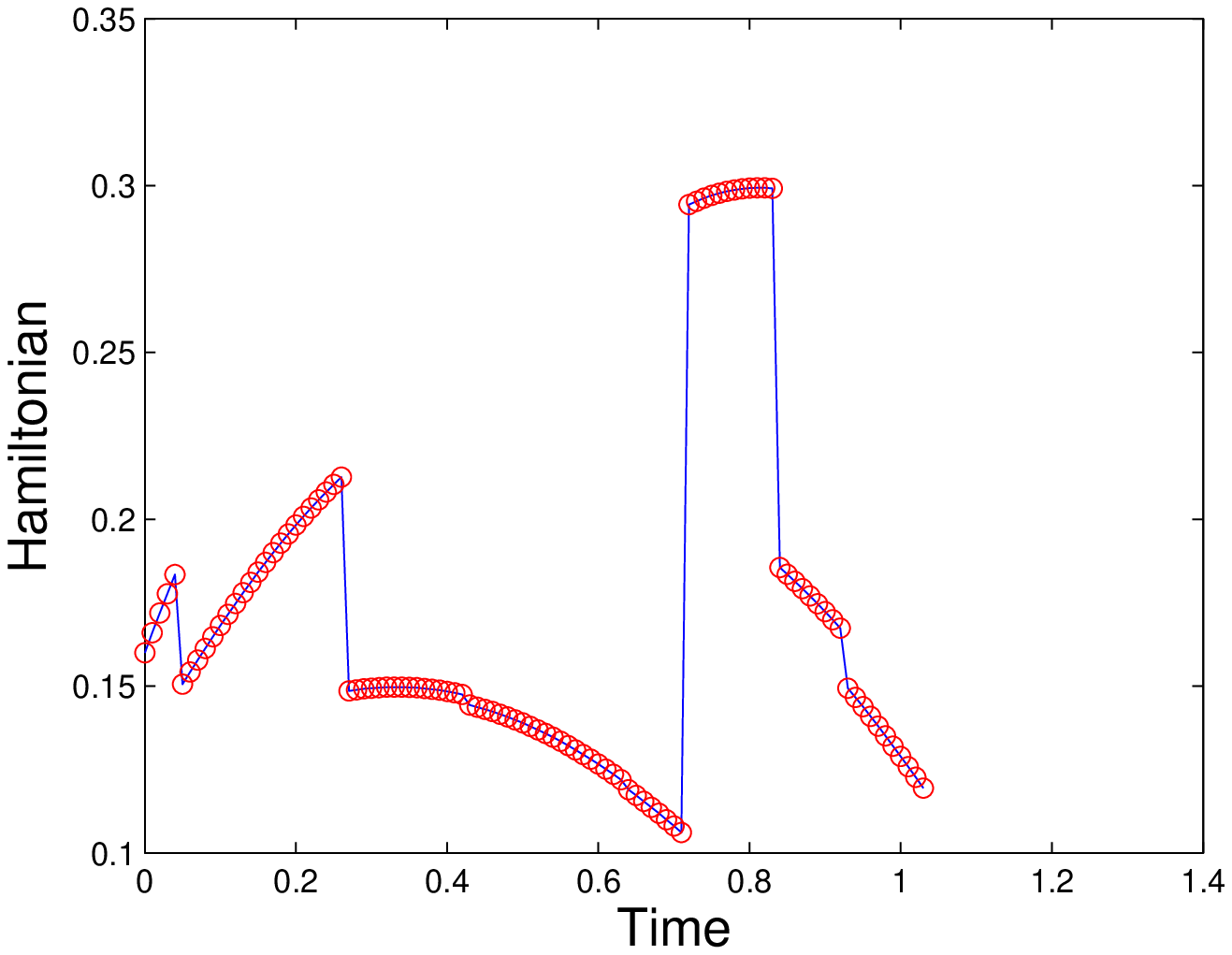}
    \end{minipage}
    \begin{minipage}{6.5cm}
       \includegraphics[width=3.8in, height=2.60in]{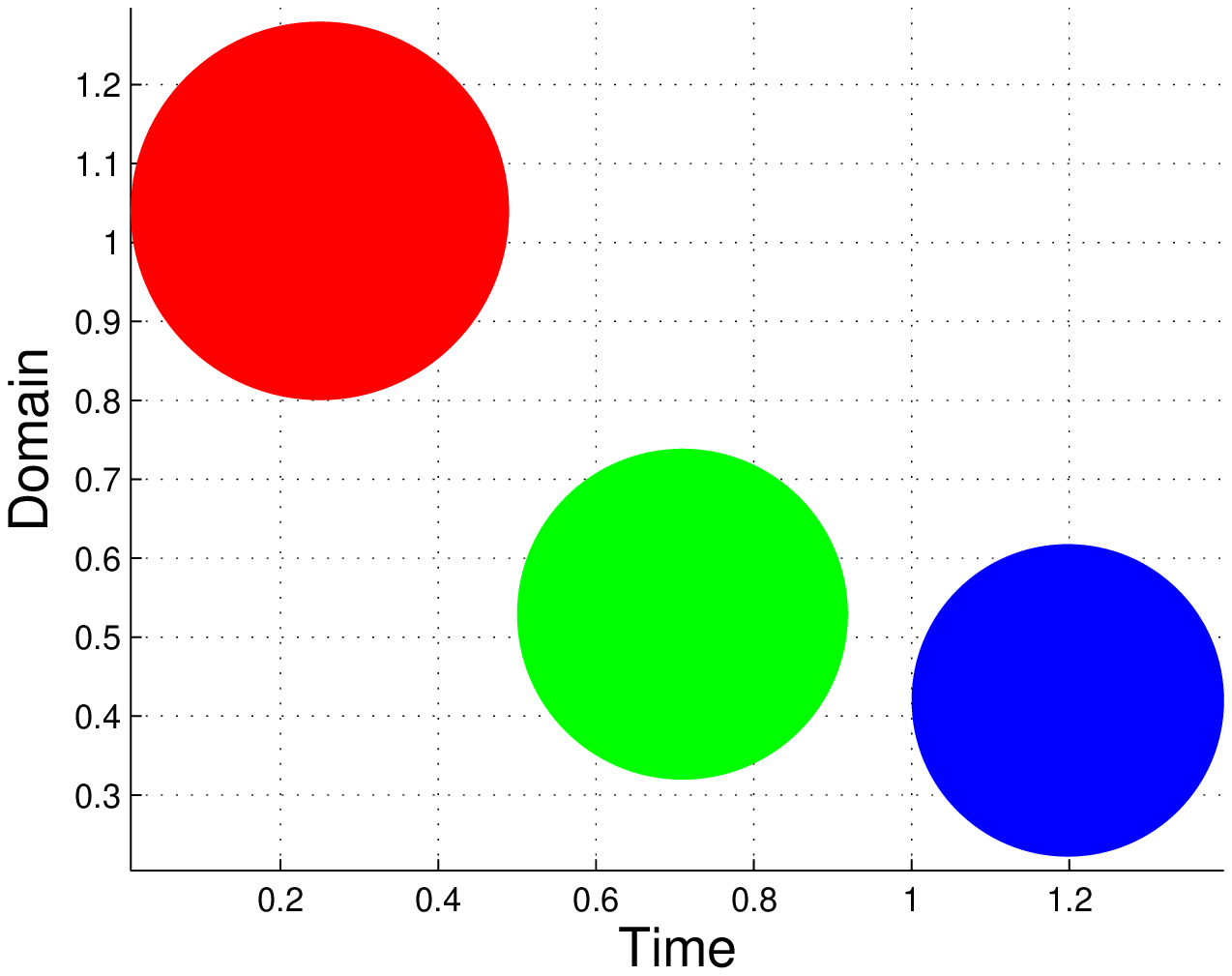}
    \end{minipage}
    \par{\scriptsize  Fig.6. Part zoom of one specific realization of the Hamiltonian and domains of the solution of SDEs $(\ref{5.1})$ by SES.}
\end{figure}

\section{Conclusion}
\ \ \ \  The main results of this paper are the construction, the convergence analysis and the numerical implementation of SES for Hamiltonian SDEs with the additive L$\acute{e}$vy noise in the Marcus form. It focuses on the mathematical approaches to preserve the symplectic structure and to realize SES. The results show that the method is effective and the numerical experiments are performed and match the results of theoretical analysis almost perfectly. More high-performance symplectic schemes, such as symplectic Runge-Kutta scheme, for Hamiltonian SDEs in the Marcus form, and numerical methods for Hamiltonian SDEs with multiplicative L$\acute{e}$vy noise in the sense of Marcus integral will be shown in our further work.

\section*{Statements}
All data in this manuscript is available. And all programs will be available on the WEB Github.

\section*{Acknowledgments}

All data in this manuscript is available. This work is supported by NSFC(No. 61841302 and 11771449). This work is also supported by the Science Research Projection in the Education Department of Fujian Province, No. JT180122, Education Reform Fund of Fujian Agriculture and Forestry University, No. 111418136. Qingyi Zhan would like to thank the Department of Applied Mathematics, Illinois Institute of Technology for the hospitality during his visit(2019-2020). And he wants to
thank Prof. Jinqiao Duan, Prof. Xiaofan Li and all members in Lab. for stochastic dynamics and computation of IIT for many fruitful discussions
during that period. Qingyi Zhan would also like to acknowledge
the sponsorship of the China Scholarship Council, CSC No. 201907870004.

\begin{backmatter}

\section*{Competing interests}
The authors declare that they have no competing interests.

\section*{Author details}
1.College  of Computer and Information Science, Fujian Agriculture and Forestry  University, Fuzhou, Fujian, 350002, PR China;

2.Department of Applied Mathematics, Illinois Institute of Technology, Chicago, IL, 60616, USA ; \\
*:Corresponding Author:Q.Zhan: zhan2017@fafu.edu.cn,qzhan3@iit.edu. Co-author: Jinqiao Duan:duan@iit.edu; Xiaofan Li: lix@iit.edu.


\bibliographystyle{bmc-mathphys} 
\bibliography{bmc_article}      


\end{backmatter}
\end{document}